\documentstyle[11pt]{article}
   \csname @addtoreset\endcsname{equation}{section}
    \csname @addtoreset\endcsname{figure}{section}
    \csname @addtoreset\endcsname{table}{section}

\newcounter{thanksnum}
\def\thanksnumber#1
{\setcounter{thanksnum}{\value{footnote}}\setcounter{footnote}{#1}%
                     \addtocounter{footnote}{-1}\footnotemark
                     \setcounter{footnote}{\value{thanksnum}}}
\def\newtheoremz#1{\@ifnextchar[{\@othmz{#1}}{\@nthmz{#1}}}

\def\@nthmz#1#2{%
\@ifnextchar[{\@xnthmz{#1}{#2}}{\@ynthmz{#1}{#2}}}

\def\@xnthmz#1#2[#3]{\expandafter\@ifdefinable\csname #1\endcsname
{\@definecounter{#1}\@addtoreset{#1}{#3}%
\expandafter\xdef\csname the#1\endcsname{\expandafter\noexpand
  \csname the#3\endcsname \@thmcountersepz \@thmcounterz{#1}}%
\global\@namedef{#1}{\@thmz{#1}{#2}}\global\@namedef{end#1}{\@endtheoremz}}}

\def\@ynthmz#1#2{\expandafter\@ifdefinable\csname #1\endcsname
{\@definecounter{#1}%
\expandafter\xdef\csname the#1\endcsname{\@thmcounterz{#1}}%
\global\@namedef{#1}{\@thm{#1}{#2}}\global\@namedef{end#1}{\@endtheoremz}}}

\def\@othmz#1[#2]#3{\expandafter\@ifdefinable\csname #1\endcsname
  {\global\@namedef{the#1}{\@nameuse{the#2}}%
\global\@namedef{#1}{\@thmz{#2}{#3}}%
\global\@namedef{end#1}{\@endtheoremz}}}

\def\@thmz#1#2{\refstepcounter
    {#1}\@ifnextchar[{\@ythmz{#1}{#2}}{\@xthmz{#1}{#2}}}

\def\@xthmz#1#2{\@begintheoremz{#2}{\csname the#1\endcsname}\ignorespaces}
\def\@ythmz#1#2[#3]{\@opargbegintheoremz{#2}{\csname
       the#1\endcsname}{#3}\ignorespaces}

\def\@thmcounterz#1{\noexpand\arabic{#1}}
\def\@thmcountersepz{.}
\def\@begintheoremz#1#2{ \trivlist \item[\hskip \labelsep{\bf #1\ #2}]}
\def\@opargbegintheoremz#1#2#3{ \trivlist
      \item[\hskip \labelsep{\bf #1\ #2\ (#3)}]}
\def\@endtheoremz{\endtrivlist}

\newtheorem{theorem}{Theorem}[section]
\newtheorem{lemma}{Lemma}[section]

\newtheorem{corollary}{Corollary}[section]
\newtheorem{condition}{Condition}[section]

\newtheorem{remark}{Remark}[section]

  \def\YY{Y}

\def\e{\varepsilon}

\def\defi{\stackrel{{\scriptscriptstyle \Delta}}{=}}

\def\N{\mu}
\def\a{\alpha}

\def\Y{{\cal Y}}

\def\w{\widehat}

\def\mes{{\rm mes\,}}
\def\arctg{{\rm arctg\,}}

\def\const{{\rm const\,}}

\def\R{{\bf R}}
\def\E{{\bf E}}
\def\P{{\bf P}}

\def\g{\gamma}
\def\C{{\bf C}}

\def\ww{\widetilde}

\def\oo{\bar}

\def\p{\partial}

\def\YY{{\w Y}}
\def\XX{{\w X}}

\def\LLL{{\cal L}}
\def\LL{{{\cal L}}}
\def\N{{\cal N}}

\newcommand{\be}{\begin{equation}}
\newcommand{\ee}{\end{equation}}
\newcommand{\baa}{\begin{eqnarray}}
\newcommand{\eaa}{\end{eqnarray}}
\newcommand{\bd}{\begin{displaymath}}
\newcommand{\ed}{\end{displaymath}}
\newcommand{\ba}{\begin{array}{ll}}
\newcommand{\ea}{\end{array}}
\newcommand{\baaa}{\begin{eqnarray*}}
\newcommand{\eaaa}{\end{eqnarray*}}\font\sm=cmr10

\date{1997 }
\title{Cordes conditions and some alternatives
for parabolic equations and discontinuous diffusion\thanks{ {\it
Differential equations } (1997) {\bf 33} (4), English translation:
pp. 433-442, in Russian: 1996, pp. 552-531.
}}
 \author{
Nikolai Dokuchaev\\
  {\sm Faculty of Mathematics and Mechanics, St.Petersburgh State University
}}
\begin{document}
\maketitle
\section*{Introduction}
We consider parabolic equations in nondivergent form with
discontinuous coefficients at higher derivatives. Their
investigation is most complicated because, in general, in the case
of discontinuous coefficients, the uniqueness of a solution for
nonlinear parabolic or elliptic equations can fail, and there is
no a priory estimate for partial derivatives of a solution. There
are some conditions that ensure regularity of solutions of
boundary value problems for second order equations and that are
known as Cordes conditions (see Cordes (1956)). These conditions
restricts the scattering  of the eigenvalues  of the matrix of the
coefficients at higher derivatives.  Related conditions from
 Talenti (1965), Koshelev (1982), Kalita (1989),  Landis (1998), on
the eigenvalues are also called Cordes type conditions.   Gihman
and Skorohod (1975) obtained  a closed
condition implicitly  as a part of the proof of the uniqueness of a weak
solution in Section 3 of Chapter 3. Cordes (1956) considered elliptic equations.   Landis
(1998) considered both elliptic and parabolic equations.  Koshelev
(1982)  considered systems of elliptic equations of divirgent type
and H\"older property of solutions.  Kalita (1989) considered
union of divergent and nondivirgent cases.
\par
Conditions from  Cordes (1956) are such that they are not
necessary satisfied even for constant non-degenerate matrices $b$,
therefore, the condition for $b=b(x)$ means   that the
corresponding inequalities are satisfied for all $x_0$ for some
non-degenerate matrix $\theta(x_0)$  and
 $\ww  b(x)=\theta (x_0)^T  b(x) \theta (x_0)$,
where    $x$  is from $\e$-neighborhood of  $x_0$  ($\e>0$ is
given). We found another condition (Condition \ref{cond3.1.A}
below) that ensures  solvability and uniqueness for first boundary
value problem for nondivirgent parabolic equation with
discontinuous diffusion coefficients. This condition ensures
existence of $ L_2$-integrable derivatives for the solution for
$L_2$-integrable free term. Prior estimate is proved, in contrast
with the existing  literature.
\par
For discontinuous diffusions, uniqueness of a weak solution cannot
be guarantied for the general case (some cases of uniqueness are
described in Gihman and Skorohod (1975),  Krylov (1980), Anulova {\it et al} (1998),
 Liptser
and Shiryaev (2000). We obtain some new conditions of uniqueness closed to
conditions Gihman and Skorohod
(1975) but sometimes less restrictive, as is shown by an example.
\subsubsection*{Some definitions} \par
Assume that we are given $T>0$ and an open  domain  $D\subset \R^n$
such that either $D=\R^n$ or $D$ is bounded with the boundary $\p D$ that  is
either  $C^2$-smooth (or such as described in  Chapter III.8 in
Ladyzhenskaya and Ural'tseva  (1968)).
\par
 We denote Euclidean norm as
$|\cdot|$, and $\bar D$ denotes the closure of a region $D$.
\par
We denote by $\|\cdot\|_{ X}$ the norm in a linear normed space
$X$, and
 $(\cdot, \cdot )_{ X}$ denotes the scalar product in  a Hilbert space $
X$.
\par
Introduce some spaces of
functions. Let $G\subset \R^k$ be an open domain, then
${W_q^m}(G)$ denotes the Sobolev  space of functions that belong
$L_q(G)$ together with first $m$ derivatives, $q\ge 1$.
\par Let $H^0\defi L_2(D)$ be the Hilbert space of complex valued functions,
and let $H^1\defi \stackrel{\scriptscriptstyle 0}{W_2^1}(D)$ be
the closure in the ${W}_1^1(D)$-norm of the set of all smooth
functions that vanish in a neighborhood of $\p D$, $k=1,2$.
Let
$H^2=W^2_2(D)\cap H^1$ be the space equipped with the norm of
$W_2^2(D)$.
\par
Let $\ell_{m}$ denotes the Lebesgue measure in  $\R ^m$, and let $
{\overline{\cal B}}_{m}$ be the $\sigma $--algebra of the Lebesgue
sets in $\R ^m$.
\par  We shall use spaces\\
$X^k=L^2([0,T],\oo{\cal B}_1,{\ell}_1,H^k)$, ${\cal
C}^k=C([0,T]; H^k)$, $k=0,1,2$, ${Y}^k=X^k\cap {\cal
C}^{k-1}$,  $k=1,2$, with the norm $\| v \|_{{Y}^k}= \| v
\|_{X^k}+ \| v \|_{{\cal C}^{k-1}}.$
\section{Solvability of boundary value problem}
Consider the domain  $D\subset\R^n$ such as described above,
$n\ge 1$. Let $Q\defi D\times [0,T]$, where
$T>0$ is given.
\par
 \par Let \be Av=\sum _{i,j=1}^n
b_{ij}(x,t)\frac{\p^2v}{\p  x_i \p x_j}  (x) + \sum_{i=1}^n
f_i(x,t)\frac{\p v}{\p x_i}(x) -\lambda (x,t)v(x), \label{3.1.1}
\ee where $(x,t)\in Q$.
\par
We are studying the problem in  $Q$ \be \left\{\ba\frac{\p   v}{\p
t}+Av=-\varphi, \\ v(x,t)| _{x\in\p D}=0, \quad   v(x,T)=\Phi
(x).\ea\right. \label{3.1.2} \ee Here   $b(x,t):\R^n\times \R\to
\R^{n\times n}$, $f(x,t):\R^n\times \R\to \R^n$, and $\lambda
(x,t):\R^n\times \R\to \C$ are measurable bounded functions,
$b_{ij}$, $f_i$, and $x_j$ are the components of $b,f$, and $x$.
\par
If $D=\R^n$, then the
boundary condition for $\p D$  vanish in (\ref{3.1.2}).
\par We assume that $b(x,t)$, $f(x,t)$, $\lambda(x,t)$ vanish for
$(x,t)\notin D\times [0,T]$.\par
\par Let us state the main conditions imposed on the matrix $b$.
\begin{condition}
\label{cond3.1.A}  The matrix  $b=b^\top$ is symmetric and has the
form $b(x,t)=\oo b(x,t)+\w b(x,t), $ where $\oo b(x,t)=\oo
b(x,t)^\top$ is a continuous bounded matrix such that
$$\delta\defi
\inf_{(x,t)\in  Q, \  \xi \in   \R^n} \frac{\xi^T\oo
b(x,t)\xi}{|\xi |^2}> 0.
$$
The matrix function $\w b(x,t)\in L^{\infty}(Q;\R^{n\times n})$ is
symmetric and such that there exists a set ${\cal
N}\subseteq\{1,\ldots,n\}$ such that
$$\w b_{ij}\equiv \w b_{ji}\equiv 0 \quad \forall i,j: i\notin
{\cal N}, j\notin  {\cal  N},
$$ and there exists a set $\{\gamma _k\}_{k\in {\cal  N}}$ such that $\g_k\in(0,2)$ for all $k$
and
$$
\!\! \w \nu= \left(\sum_{k\in\N}\frac{1}{2\gamma_k}\right)
\mathop{{\rm ess\, sup}}_{x,t} \sum_{k\in \N} \Biggl(\sum_{i\in\N}
\w b_{ik}(x,t)^2  +\!4\sum_{i\notin\N}\w b_{ik}(x,t)^2
\!+\!\frac{\gamma_k}{2-\gamma_k}\w b_{kk}(x,t)^2\Biggr) <\!\delta
^2.
$$
\end{condition}
\begin{remark} If ${\rm card}\,\N<n$, then Condition \ref{cond3.1.A}  allows
bigger than for $\N=\{1,\ldots.,n\}$  values $\w b_{ij}$ for $i\in
\N$, $j\in \N$. Different $\gamma_k$ also make this condition less
restrictive: for instance, if $\w b_{kk}\equiv 0$, then we can
allow   $\gamma _k=2-0$. \end{remark} In particular, the condition
for $\w \nu$ is satisfied if
$$
\mathop{{\rm ess\, sup}}_{x,t} \sum^n_{i,k=1}\w b_{ik}(x,t)^2
<\frac{\delta ^2}{n}.
$$\par
The next condition is not so principal, since it  deals with low
order coefficients and the continuous part $\oo b$.
\begin{condition}\label{cond3.1.B}
There exists a domain $D_1\subseteq \R^n$ and functions
$b^{(\e)}(x,t):\R^n  \times \R\to \R^{n \times   n}$,
$f^{(\e)}(x,t):\R^n \times \R\to \R^n$, $\lambda ^{(\e)}(x,t):\R^n
\times \R\to \C$, $\e >0$,  such that $\mes  D_1<+\infty$,
\begin{eqnarray}
\nu_b  (\e)= \|b^{(\e)}-\oo b\|_{L_{\infty}(Q)} \to 0 \quad
\mbox{as} \quad \e\to 0, \quad
\nonumber
\\
\oo \nu_b(\e)= \mathop{\rm ess\, sup}_{(x,t)\in    Q} | \frac{\p
b^{(\e)}}{\p x}(x,t)|  <+\infty \quad \forall \e >0, \nonumber\\
\nu _f (\e)=\|f^{(\e)}-f\| _{L_n(Q_1)} +\mathop{\rm ess\,
sup}_{(x,t)\in  Q\backslash Q_1} |f^{(\e)}(x,t)-f(x,t)| \to 0
\quad \mbox{ as} \quad  \e\to 0, \quad \nonumber\\ \oo \nu_f(\e )=
\mathop{\rm ess\, sup}_{(x,t)\in    Q} | \frac{\p f^{(\e)}}{\p
x}(x,t)| <+\infty \quad \forall  \e >0, \nonumber\\
\nu_\lambda(\e)=\|\lambda^{(\e)}-\lambda\|_{L_r(Q_1)} +
\mathop{\rm ess\, sup}_{(x,t)\in   Q\backslash Q_1}
|\lambda^{(\e)}(x,t)-\lambda(x,t)| \to 0 \quad \mbox{ as} \quad
\e\to 0, \quad \nonumber\\ \oo \nu_\lambda(\e   )   = \mathop{\rm
ess\, sup}_{(x,t)\in    Q} | \frac{\p \lambda^{(\e)}}{\p x}(x,t)|
<+\infty \quad \forall \e
>0
\nonumber\end{eqnarray} Here $Q_1\defi D_1\times (0,T)$, $r\defi
\max (1,n/2)$.
\end{condition}
\par
\begin{remark} Condition  \ref{cond3.1.B} is satisfied if  $f,\oo b,
\lambda$ are bounded and  $D$ is bounded. In that case,  we can
take  $D_1=D$ and the Sobolev averages of the functions $\oo
b,f,\lambda $ as $b^{(\e)},f^{(\e)},\lambda^{(\e)} $ respectively.
Note that Condition  \ref{cond3.1.B} implies that
$$
\|f\| _{L_n(Q_1)}
+\mathop{\rm ess\, sup}_{(x,t)\in  Q\backslash Q_1}
|f(x,t)|<+\infty,
\quad
\quad\|\lambda\|_{L_r(Q_1)}
+
\mathop{\rm ess\, sup}_{(x,t)\in   Q\backslash Q_1}
|\lambda(x,t)|<+\infty .
$$
\end{remark}
\par
We introduce the set of parameters
$$
\ba {\cal P} \defi \biggl( n,\,\, D,\,\, T,\,\, \delta ,\,\, {\cal
N},\,\, \{\gamma_k\}_{k\in \N},\,\,\\ \sup_{x,t} |{b(x,t)}|,\,\,
\sup_{x,t} |{f(x,t)}| ,\,\, \sup_{x,t} | \lambda (x,t|,\,\, \w \nu
,\,\, \nu _b(\cdot ),\,\, \oo\nu_b(\cdot ), \,\, \nu _f(\cdot
),\,\, \oo \nu_f(\cdot ),\,\, \nu_{\lambda}(\cdot ),\,\, \oo
\nu_{\lambda}(\cdot ) \biggr). \ea
$$ We have that ${\cal P} $ includes $\nu_b(\cdot)$,
hence  ${\cal P}$ depends on the modulus of continuity of $\oo b$.
\begin{theorem}\label{Th3.1.1}
 Assume that Conditions \ref{cond3.1.A}--\ref{cond3.1.B} are satisfied.
 Then problem (\ref{3.1.2}) has
the unique solution  $v\in {\Y}^2$ for any $\varphi \in L_2(Q)$,
$\Phi \in H^1$, and
\be
\| v \|_{{\Y}^2}
 \le   c(  \|\varphi \|_{L_2(Q)}  +
\|\Phi \|_{H^1}), \label{3.1.3} \ee where  $c=c({\cal P})$ is a
constant that depends on ${\cal P} $.
\end{theorem}
 \par
We shall need some auxiliary spaces to prove the theorem. Let $\w
H^2$ be the set of  $v\in W_2^2(D)\cap  H^1$ with the special norm
\be \| v \|_{\w H^2}= \left(\sum_{k\in \N} \left\{ \sum_{i=1}^n
\left\| \frac{\p ^2v}{\p x _k \p x_i}\right\|_{H^0}^2-
\frac{\gamma_k}{2} \left\|\frac{\p ^2v}{\p x _k^2}
\right\|_{H^0}^2\right\}\right)^{1/2}+\alpha_1 \| v\|_{W^2_2(D)}.
\label{3.1.4} \ee Here $\alpha _1>0$ is some constant. \par
Introduce Banach spaces $\XX^2=L^2([0,T],\oo{\cal B}_1,{\ell}_1,
\w H^2)$ and ${\YY}^2=\XX^2\cap {\cal C}^1$ with the norm \be \| v
\|_{{\YY}^2}= \| v \|_{\XX^2} + \alpha _2  \| v \|_{{\cal C}^1}.
\label{3.1.5} \ee Here $\alpha _2>0$ is a constant.
\begin{remark}\label{rem3.1.3} Since $\g_k\in(0,2)$ for all $k$,
({3.1.4}) defines a norm, the norm $\w H^2$ is equivalent to
the norm
 $W^2_2(D)$, and the norm $\YY^2$ is equivalent to the norm  ${Y}^2$.
\end{remark}
Therefore, to prove Theorem \ref{Th3.1.1}, it suffices to  prove
the following   theorem.
\begin{theorem}
\label{Th3.1.2} Assume that Conditions
\ref{cond3.1.A}-\ref{cond3.1.B} are satisfied. Then problem
(\ref{3.1.2}) has an unique solution
 $v\in {\YY}^2$ for any
$\varphi \in L_2(Q)$ and $\Phi \in H^1$, and \be \| v \|_{{\YY}^2}
\le   c(  \|\varphi \|_{L_2(Q)}  + \|\Phi \|_{H^1}), \label{3.1.6}
\ee where  $c>0$ is a constant that depends only on ${\cal P} $
and $\alpha_1$, $\alpha_2$.
\end{theorem}
\begin{remark} For $D=\R^n$
a closed to Theorem  \ref{Th3.1.1} was announced   in Dokuchaev
(1996), where, however, the estimate was obtained for the
derivatives with discontinuous coefficients only, just to make the
equation meaningful).
\end{remark}
\section{Examples}
Let $b=b(x)$, and let $\lambda_1$,\ldots,$\lambda_n$ be its
eigenvalues. The classic Cordes conditions from Cordes (1956) was
formulated for $n\ge 3$ as \be \exists  \e   >0: \quad
(n-1)\sum_{i<j}\left(\lambda_i-
     \lambda_j\right)^2<(1-\e )\left(\sum_{i=1}^n\lambda_i\right)^2.
\label{3.2.1} \ee It was shown by Talenti (1965)   that
(\ref{3.2.1}) is equivalent to \be \exists \e   >0:     \quad
(n-1+\e)\sum_{i=1}^n\lambda_i ^2 = (n-1+\e)\sum_{i,j=1}^n b_{ij}^2
<\left(\sum_{i=1}^n
b_{ii}\right)^2=\left(\sum_{i=1}^n\lambda_i\right)^2.
\label{3.2.1'} \ee This form (\ref{3.2.1'}) can be given also to
the condition from Kalita (1989) for a system with one
nondivirgent equation.
\par
Conditions from Landis (1998)  has the form \be \exists \e >0:
\quad \sum_{i=1}^n \lambda_i< (n+2-\e)\min
\left\{\lambda_1,\ldots,\lambda_n\right\}. \label{3.2.2} \ee The
condition  from Section 3, Chapter 3 from Gihman and Skorohod
(1975) is such that in the simplest case can be written as \be
\exists \e
>0:     \quad {\rm Tr\,}\bigl( (b-I)^2\bigr)<1-\e. \label{3.2.3}
\ee (In Gihman and Skorohod (1975), $I$ was replaced   for a
smooth matrix function).
 \par
In our notations, the last condition can be rewritten as \be \oo
b\equiv I, \quad \exists  \e >0: \quad \sum_{i,j=1}^n\w
b_{ij}^2<1-\e. \label{3.2.3'} \ee The regularity of the parabolic
equation established by Gihman and Skorohod (1975) under condition
(\ref{3.2.3}) is weaker than the regularity established by Theorem
\ref{Th3.1.1} \par Note that Gihman and Skorohod (1975) obtained
the regularity that was just enough to ensure the uniqueness of a
weak solution of some Ito's equation. In fact, conditions
(\ref{3.2.3}), (\ref{3.2.3'}) are sufficient for Theorem
\ref{Th3.1.1} as well. We leave it without proof; note that
there is a proof similar to the proof given below and different from the one given in
Gihman and Skorohod
(1975).
\par
In fact, Cordes conditions mean that inequalities
  (\ref{3.2.1})--(\ref{3.2.2})
are satisfied for all $x_0$ for some non-degenerate matrix $\theta
(x_0)$ and for all matrices $\w  b(x)=\theta  (x_0)^T b(x)\theta
(x_0)$  ,  where   $x$ is from the   $\e$-neighborhood of $x_0$,
and where $\e >0$ is given. Similarly,  condition (\ref{3.2.2})
was adjusted in Landis (1998), and condition (\ref{3.2.3}) was
adjusted in Gihman and Skorohod (1975).
\par
Let $n=3$, $b(x,t)\equiv b(x)$,
     $$b(x)=\left(\begin{array}{ccc} 1&\alpha(x)&\beta(x)\cr
     \alpha(x)&1&0\cr
      \beta(x)&0&1 \end{array} \right), \quad
     \widehat b(x)=\left(\begin{array}{ccc}
     0&\alpha(x)&\beta(x)\cr   \alpha(x)&0&0\cr
     \beta(x)&0&0 \end{array} \right), $$
where  $\alpha(x)$,$\beta(x)$ are arbitrary measurable functions,
$|\alpha(x)|\le\alpha=\const$,   $|  \beta(x)|    \le
\beta=\const$,   and functions    $\alpha(x)$,$\beta(x)$ are quite
irregular.
\par
It is easy to see that Condition  \ref{cond3.1.A}
 is satisfied if $\alpha^2+\beta^2<1$ for  $\N=\{1\}$ and for some
$\gamma_1<2$ being close enough to 2.
\par
The spectrum of  $b$ is $\{1$, $1-\sqrt{\alpha(x)^2+\beta(x)^2}$,
$1+\sqrt{\alpha(x)^2+\beta(x)^2}\}$. Then conditions (\ref{3.2.1}),
(\ref{3.2.1'}) fails if $(\oo\alpha^2+\oo\beta^2)\ge  3/4$,
and (\ref{3.2.2}) fails if   $(\oo\alpha^2+\oo\beta^2)\ge
2/5$. Conditions (\ref{3.2.3}) and (\ref{3.2.3'}) fail if $\oo
\alpha^2+\oo \beta^2>1/2$.
\par
Therefore,  Condition \ref{cond3.1.A} is less restrictive for this
example than condition (\ref{3.2.3'}) or the conditions from
Cordes (1956), Gihman and Skorohod (1975), Kalita (1989),
Koshelev (1982), Landis (1998), Talenti (1965).
\par
There may be opposite  examples when condition (\ref{3.2.1}) is
satisfied, but Condition \ref{cond3.1.A}  fails.
\section{Proof of Theorem 3.1.2.}
The main idea is  to prove theorem  for some $\e=\e({\cal P})>0$
for $u$ replaced with \be u_{\e}(x,t)\defi u(x,t)\exp\{K(\e)t\},
\label{3.3.1} \ee where   $K(\e)>0$ is a function of   $\e$ such
that  $K(\e)  \to  +\infty $ as  $\e  \to  +0$.
\par
Let $\ww \lambda^{(\e)}(x,t)\defi \lambda^{(\e)}(x,t)+K(\e)$, and let
$$
A_{\e}u\defi \sum_{i,j=1}^{n}b^{(\e)}_{ij}(x,t) \frac{\p^2 u}{\p
x_i\p x_j}(x)+ \sum_{i}^{n} f^{(\e)}_{i}(x,t) \frac{\p u}{\p
x_i}(x)-\ww\lambda^{(\e)}(x,t)u(x,t).
$$
Consider the problem \be\left\{\ba\frac{\p v}{\p t}+
A_{\e}v=-\varphi , \\ v(x,t)|_{x\in \p D}=0, \quad v(x,T)=\Phi
(x).\ea\right. \label{3.3.2} \ee Introduce the operators
$L(\e):X^0 \to {\YY}^2$, $\LL(\e):H^1  \to {\YY}^2$ such that
$u \defi L(\e)\varphi+\LL(\e)\Phi $ is the solution of (\ref{3.3.2}).
Let $\|L(\e)\|$ denotes the norm of the operator $L(\e):X^0 \to
{\YY}^2$, and let $\|\LL(\e)\|$ denotes the norm of the operator
$\LL(\e):H^1 \to {\YY}^2$.
\begin{lemma}\label{lemma3.3.1}
For any $ \gamma>0$, there exists   a small enough   $\e_*>0$, and
a function $K(\e)>0$ (increasing as  $\e \to 0$), and
$\alpha_i=\alpha_i(\gamma,{\cal P})$, $i=1,2$, in
(\ref{3.1.4})-(\ref{3.1.5}), such that
 $\e_*=\e_*(\gamma,{\cal
P})$, $K(\cdot)=K(\cdot,\gamma,{\cal P})$, and \be\|L(\e)\|\le
\gamma+ \frac{1}{\delta}\biggl(
\sum_{k\in\N}\frac{1}{2\gamma_k}\biggr)^{1/2},
 \quad  \|\LL(\e)\|\le c_0  \quad \forall \e \in (0,\e_*],
\label{3.3.3} \ee where  $c_0=c_0({\cal P},\a_1,\a_2)$ is a
constant.
\end{lemma}
\par
{\it Proof}.   Let  $\varphi \in X^0$ be a smooth function with a
compact support inside  $Q$. Set $\g_k\defi 1$ for $k\notin \N$.
 \par
Let $v=L(\e)\varphi$. We have \baaa
\frac{1}{2}\|v(\cdot,t_1)\|^2_{H^0}=\frac{1}{2}\|v(\cdot,T)\|^2_{H^0}
+ \int_{t_1}^T (v,A_{\e}v+\varphi)_{H^0} ds. \label{3.3.4)} \eaaa We
shall use below  the obvious inequality
$$
2\alpha\beta  \le    \e\alpha^2+\e^{-1}\beta   ^2 \quad \forall
\alpha,\beta,\e \in \R, \e>0.
$$
In particular, \baaa ( v, \varphi)_{H^0} \le \frac{1}{2\e _1}
\|v\|^2_{H^0}+\frac{\e_1}{2}\|\varphi\|^2_{H^0}\quad \forall
\e_1>0. \label{3.3.5} \eaaa We have the  estimate
\begin{eqnarray}
(v,A_{\e}v+\varphi)_{H^0} = \biggl(
v,\sum_{i,j=1}^{n}b^{(\e)}_{ij}(\cdot,t) \frac{\p^2 v}{\p x_i\p
x_j} +\sum_{i=1}^{n} f_i^{(\e)}(\cdot,t)\frac{\p v}{\p x_i}
 -\ww \lambda^{(\e)}(\cdot,t)v(\cdot,t) \biggr)_{H^0}
 \nonumber\\
=\sum_{i,j=1}^{n}\biggl\{-\biggl(v,\frac{\p b_{ij}^{(\e)}}{\p x_j}
\frac{\p v}{\p x_i} \biggr)_{H^0} -\biggl(\frac{\p v}{\p
x_j},b^{(\e)}_{ij} \frac{\p v}{\p x_i}\biggr)_{H^0}\biggr\}
\nonumber\\
- \frac{1}{2}\biggl( v^2,\sum_{i=1}^{n}
\frac{\p f_i^{(\e)}}{\p x_i}\biggr)_{H^0}
-(v,\lambda^{(\e)}v)_{H^0} -K(\e)\|v\|^2_{H^0} +
(v,\varphi)_{H^0}
\nonumber\\
\le (-\delta+\nu_1) \sum_{j=1}^{n} \biggl\|\frac{\p v}{\p
x_j}\biggr\|^2_{H^0} -K(\e)\|v\|^2_{H^0} +c_1 \|v\|^2_{H^0}+
\frac{\e_1}{2}\|\varphi\|^2_{H^0}, \nonumber\end{eqnarray} where
$\e_1>0$, $\nu_1>0$ can be arbitrarily small, and $c_1$ depends on
$\e$, $\e_1$, $\nu_1$,  ${\cal P}$.  Hence we have that choosing
$K(\e)=K(\e,\nu)>c_1$ for $\nu>0$ can ensure that \baa
 \|L(\e)\varphi\|_{{\YY}^1}   \le   \nu\|\varphi\|_{X^0}    \quad
\forall  \e \in (0,\e_*],\ \  \forall \varphi \in X^0.
\label{3.3.6} \eaa We have that \baaa \left\|\frac{\p v}{\p
x_k}(\cdot,t_1)\right\|^2_{H^0} -\left\|\frac{\p v}{\p
x_k}(\cdot,T)\right\|^2_{H^0} = 2\int_{t_1}^T \left(\frac{\p v}{\p
x_k}, \frac{\p}{\p x_k} \bigl(A_{\e}v+\varphi\bigr)\right)_{H^0}
ds. \label{3.3.7} \eaaa Remind that $\varphi$ has compact support
inside $Q$. Then \baaa \biggl( \frac{\p v}{\p x_k}, \frac{\p
\varphi}{\p x_k}\biggr)_{H^0} \le \frac{\delta \gamma_k}{2}
\biggl\| \frac{\p^2 v}{\p x_k^2}\biggr\|^2_{H^0} +
\frac{1}{2\delta \gamma_k}\|\varphi\|^2_{H^0} . \label{3.3.8} \eaaa
\par
Note that if  $b^{(\e)}\in C^2$ then \be \ba &\biggl( \frac{\p
v}{\p x_k}, \frac{\p b^{(\e)}_{ij}}{\p x_k} \frac{\p^2 v}{\p x_i\p
x_j}\biggr)_{H^0}
\\&=-\biggl( \frac{\p v}{\p x_k}, \frac{\p^2
b^{(\e)}_{ij}}{\p x_i\p x_k} \frac{\p v}{\p x_j}\biggr)_{H^0}
-\biggl( \frac{\p^2 v}{\p x_i\p x_k}, \frac{\p b^{(\e)}_{ij}}{\p
x_k} \frac{\p v}{\p x_j}\biggr)_{H^0}+\int_{\p D}\w J_{ijk} ds \\
&\hphantom{=}= \biggl( \frac{\p v}{\p x_k}, \frac{\p
b^{(\e)}_{ij}}{\p x_i} \frac{\p^2 v}{\p x_j\p x_k}\biggr)_{H^0}
+\biggl( \frac{\p^2 v}{\p x_k^2}, \frac{\p b^{(\e)}_{ij}}{\p x_i}
 \frac{\p v}{\p x_j}\biggr)_{H^0} -\biggl( \frac{\p^2 v}{\p x_i\p x_k},
\frac{\p b^{(\e)}_{ij}}{\p x_k} \frac{\p v}{\p x_k} \biggr)_{H^0}
+\int_{\p D}J'_{ijk}ds , \ea
\label{3.3.9} \ee where
$$
J'_{ijk}=\w J_{ijk} -
\frac{\p \oo v}{\p x_k}
\frac{\p b^{(\e)}_{ij}}{\p x_i}\frac{\p v}{\p x_j}\cos({\bf n},e_k),
\quad
\w J_{ijk} = \frac{\p \oo v}{\p x_k} \frac{\p b^{(\e)}_{ij}}{\p x_k}
\frac{\p v}{\p x_j}\cos({\bf n},e_i),
$$
${\bf n}={\bf n}(s)$ is the outward pointing normal to the surface
$\p D$  at the point $s \in \p D$, and $e_k$ is the $k$th  basis
vector in the Euclidean space $\R^n=\{ x_1,\ldots,x_n \}$.
\par
If  $b^{(\e)}$ is general, then the right hand and the left hand
expressions in (\ref{3.3.9}) are still equal. Hence, we obtain
\baaa &&\biggl( \frac{\p v}{\p x_k}, \frac{\p b^{(\e)}_{ij}}{\p
x_k} \frac{\p^2 v}{\p x_i\p x_j}\biggr)_{H^0}\\ &&\le
\e_2\biggl(\biggl\|\frac{\p^2 v}{\p x_k\p x_j}\biggr\|^2_{H^0}
+\biggl\|\frac{\p^2 v}{\p x_k\p x_i}\biggr\|^2_{H^0}
     + \biggl\|\frac{\p^2 v}{\p x_k^2}\biggr\|^2_{H^0} \biggr)
+c_2 \frac{1}{2\e_2}\|v\|^2_{H^1} +\int_{\p D}J'_{ijk} ds \quad
\forall  \e_2>0, \label{3.3.10} \eaaa where the constant $c_2$
depends only on ${\cal P}$.   \par Therefore, \baa
\biggl( \frac{\p v}{\p x_k} , \frac{\p }{\p x_k}
(A_{\e}v+\varphi)\biggr)_{H^0}\hphantom{xxxxxxxxxxxxxxxxxxxxxxxxxxxxxxxxxxxxxxxxxxxxx}
\nonumber\\
= \biggl( \frac{\p v}{\p x_k}, \frac{\p }{\p
x_k}\biggl\{\sum_{i,j=1}^n b^{(\e)}_{ij}(\cdot,t) \frac{\p^2 v}{\p
x_i\p x_j} +\sum_{i=1}^{n} f_i^{(\e)}(\cdot,t)\frac{\p v}{\p x_i}
-\ww \lambda
^{(\e)}(\cdot,t)v(\cdot,t)+\varphi(\cdot,t)\biggr\}\biggr)_{H^0}
\nonumber\\
= \sum_{i,j=1}^{n}\biggl\{ \biggl( \frac{\p v}{\p x_k},
\frac{\p b^{(\e)}_{ij}}{\p x_k}
 \frac{\p^2 v}{\p x_i\p x_j}\biggr)_{H^0}
-\biggl( \frac{\p^2 v}{\p x_k\p x_i}, b_{ij}^{(\e)}\frac{\p^2
v}{\p x_k\p x_j}\biggr)_{H^0}\biggr\}\nonumber\\
+\sum_{i=1}^{n} \biggl\{ \biggl( \frac{\p v}{\p x_k}, \frac{\p
f^{(\e)}}{\p x_k} \frac{\p v}{\p x_i}\biggr)_{H^0} + \biggl(
\frac{\p v}{\p x_k}, f_i^{(\e)} \frac{\p^2 v}{\p x_k\p x_i}
\biggr)_{H^0}\biggr\}
\nonumber\\
-\biggl( \frac{\p v}{\p x_k},
\frac{\p  \lambda^{(\e)}}{\p x_k}v+\lambda^{(\e)}
\frac{\p v}{\p x_k}\biggr)_{H^0}
-K(\e)\biggl\|\frac{\p v}{\p x_k}\biggr\|^2_{H^0}
+\biggl( \frac{\p v}{\p x_k}, \frac{\p \varphi}{\p x_k}\biggr)_{H^0}
+\int_{\p D}J_{ijk}ds
\nonumber\\
\le (-\delta+\nu_2+2\e_3)\sum_{j=1}^{n}
\biggl\|\frac{\p^2   v}{\p   x_k\p x_j}\biggr\|^2_{H^0} +
\left(\frac{\delta\gamma_k}{2} +\e_3\right)
\biggl\|\frac{\p^2 v}{\p x_k^2}\biggl\|^2_{H^0}      +
c_2\|v\|^2_{H^1}
\nonumber\\
+\sum_{i,j=1}^{n} \int_{\p D}J_{ijk} ds +
\left(\frac{1}{2\delta\g_k}+ \frac{\e_1
}{2}\right)\|\varphi\|^2_{H^0} -K(\e) \left\|\frac{\p v}{\p
x_k}\right\|^2_{H^0},\label{3.3.11} \eaa where the
constant $c_2$ depends only on   ${\cal P}$, constants $\e_3>0$
and $\nu_2>0$ can be arbitrarily small,
$$
J_{ijk} =J'_{ijk} + J''_{ijk}, \quad
J''_{ijk}=\frac{\p \oo v}{\p x_k} b^{(\e)}_{ij}
\frac{\p^2 v}{\p x_i\p x_j} \cos({\bf n},e_k).
$$
Let us estimate $\int_{\p D}J_{ijk}$. It vanishes if  $D=\R^n$ (as
well as all integrals over the boundary $\p D$).  For a bounded
domain $D$, we mainly follow the approach from  Section 3.8
Ladyzhenskaya and Ural'tseva  (1968). Let $x^0=\{x_i^0\}^n_{i=1}
\in \p D$ be an arbitrary point. In its neighborhood, we introduce
local Cartesian coordinates $y_m=\sum_{k=1}^{n}c_{mk}(x_k-x^0_k)$
such that the axis  $y_n$ is directed along the outward normal
${\bf n}={\bf n}(x_0)$ and $\{c_{mk}\}$ is an orthogonal matrix.
\par
Let $y_n=\omega (y_1,\ldots,y_{n-1})$ be an equation determining the
surface  $\p D$ in a neighborhood of the origin. By the properties
of the surface $\p D$, the first order and second order
derivatives of the function  $\omega $ are bounded. Since
$\{c_{mk}\}$ is an orthogonal matrix, we have
$x_k-x^0_k=\sum_{m=1}^{n}c_{km}y_m$.  Therefore, $\cos({\bf
n},e_m)=c_{nm}$, $m=1,\ldots,n$. Then \baaa J'_{ijk}
=\sum_{m=1}^{n}c_{mk} \frac{\p \oo v}{\p y_m}
 \sum_{p=1}^{n}c_{pi}
\frac{\p v}{\p y_p}\biggl(\sum_{q=1}^{n} \frac{\p
b^{(\e)}_{ij}}{\p y_q}c_{qk}c_{ni}-\sum_{r=1}^{n} \frac{\p
b^{(\e)}_{ij}}{\p y_r} c_{rk}c_{nk}\biggr), \nonumber\\ J''_{ijk}
=\sum_{m=1}^{n}c_{mk}\frac{\p \oo v}{\p y_m} b_{ij} c_{nk}
\sum_{p,q=1}^{n}c_{pi}c_{qj} \frac{\p^2 v}{\p y_p\p y_q} . \eaaa
The boundary condition $v(x,t)|_{x\in \p D}=0$ has the form $$
     v(y_1,\ldots,y_{n-1},\omega (y_1,\ldots,y^{n-1}),t)=0
$$
identically  with respect to $y_1,\ldots,y_{n-1}$ near the point
$y_1=\ldots=y_{n-1}=0$. Let us differentiate this identity with
respect to  $y_p$ and $y_q$, $p,q=1,\ldots,n-1$, and take into
account that $$ \frac{\p \omega}{\p y_p}=0 \quad (p=1,\ldots,n-1).
$$ at $x_0$. Then
$$
\frac{\p v}{\p y_p}=0 , \quad
\frac{\p^2 v}{\p y_p\p y_q}  =- \frac{\p v}{\p y_n}
\frac{\p^2 \omega}{\p y_p\p y_q} =
-\frac{\p v}{\p {\bf n}} \frac{\p^2 \omega}{\p y_p\p y_q}
\quad (p,q=1,\ldots,n-1).
$$
Hence \baa \int_{\p D}J_{ijk}ds \le  \w  c_1\int_{\p  D}
\biggl| \frac{\p  v}{\p  {\bf n}}\biggr|^2 ds
\le \e_4\sum_{i,j=1}^n \int_{ D} \biggl|\frac{\p^2 v}{\p x_i\p
x_j}(x)\biggr|^2dx +\w c_2(1+\e^{-1}_4) \|v\|^2_{H^1} \quad
\forall \e_4>0 \label{3.3.12} \eaa  for some constants $\w
c_i=\w c_i(\e,{\cal P})$. The last estimate follows from the
estimate (2.38) in Chapter 2 from Ladyzhenskaya and Ural'tseva
(1968).
\par
As mentioned above,  for a suitable choice of the functions
$K(\e)=K(\e,\nu)$ and for an arbitrarily small   $\nu>0$, one can
provide the estimate $\|L(\e)\varphi\|_{{\Y}^1}   \le
\nu\|\varphi\|_{X^0}$  $(\forall  \e  \in   (0,\e_*],   \forall
\varphi \in X^0)$. The constants  $\e_3>0$, $\e_4>0$, and
$\nu_2>0$ can be arbitrarily small, and the constant $c_1$ depends
on $\e , \e_1, \nu_1 ,  \gamma_k$ and , ${\cal P}$. Combining
(\ref{3.3.6}) with (\ref{3.3.11}) and (\ref{3.3.12}), we see that
for some function $K(\e)$ we have \baa \sum_{k\in\N}\biggl(
\int_{0}^{T} dt\int_{ D}\! \biggl\{
\left(\delta-\nu_2-2\e_3\right) \sum_{i=1}^{n}\left|\frac{\p^2
v}{\p x_k\p x_i}(x,t)\right|^2 -\left(\frac{\delta\gamma_k}{2}\!
+\!\e_3\!\right) \left|\frac{\p^2 v}{\p x_k^2}(x,t)\right|^2
\biggr\} dx\nonumber
\\ + \frac{1}{2}\sup_t
\left\|\frac{\p  v}{\p x_k}(\cdot,t)\right\|^2_{H^0}
\biggr)\nonumber\\     \le \sum_{k\in\N}\biggl( c_2 \|v\|^2_{X^1}
+\left( \frac{1}{2\delta\gamma_k}+
\frac{\e_1}{2}\right)\|\varphi\|^2_{X^0}\biggr) \le
\sum_{k\in\N}\left(\nu c_2+ \frac{1}{2\delta\gamma_k}+
\frac{\e_1}{2}\right)\|\varphi\|^2_{X^0}. \label{3.3.13} \eaa
Therefore, \baa \sum_{k\in\N}\biggl( \int_{0}^T\int_{ D}dt
\biggl\{ \left(\delta-\e_5\right)\sum_{i=1}^{n} \left|\frac{\p^2
v}{\p x_k\p x_i}(x,t)\right|^2 -
\frac{\delta\gamma_k}{2}\left|\frac{\p^2 v}{\p
x_k^2}(x,t)\right|^2 \biggr\}dx\nonumber\\ +\frac{1}{2}\sup_t
\left\|\frac{\p v}{\p x_k} (\cdot,t)\right\|^2_{H^0}\biggr) \le
\sum_{k\in\N} \left(\frac{1}{2\delta\gamma_k} +\e_6\right)
\|\varphi\|^2_{X^0}\nonumber \eaa for some sufficiently small
$\e_i=\e_i(\e,{\cal P})>0$, $i=5,6$. (Here $\nu_2,\e_3$ are from
(\ref{3.3.11})). Take the sum in (\ref{3.3.13}) with respect to
$k=1,\ldots,n$ and choose a sufficiently small number
$\alpha_1=\alpha_1(\g,{\cal P})$.  This, together with
(\ref{3.3.13}), yields the first estimate in (\ref{3.3.3}).
\par
In a similar way, taking into account the initial condition in
(\ref{3.3.13}) and taking the sum in (\ref{3.3.13}) with respect
to $k=1,\ldots,n$,  we obtain the estimate $ \|v\|_{{\YY}^2}\le
\ww c\|\Phi \|_{H^1}
$
for $v=\LL(\e)\Phi $, where  $\ww c=\ww c({\cal P})$ is a
constant. Then we obtain the assertion  of Lemma
\ref{lemma3.3.1}. $\Box$
\par
Introduce the operator  $R(\e): {\YY}^2\to {\YY}^2$ \baaa
R(\e)v=L(\e) \biggl\{ \sum_{i,j=1}^{n}\w b_{ij}\frac{\p^2 v}{\p
x_i\p x_j}+\sum_{i,j=1}^{n} \bigl[ \oo b_{ij}-b^{(\e)}_{ij}\bigr]
\frac{\p^2 v}{\p x_i\p x_j} +\sum_{i=1}^n \bigl[
f_i-f_i^{(\e)}\bigr] \frac{\p v}{\p
x_i}-\bigl[\lambda-\lambda^{(\e)}\bigr]v \biggr\}. \label{3.3.14}
\eaaa \begin{lemma}\label{lemma3.3.2}
 There exists a number   $\oo\e=\oo\e({\cal P})>0$ such
 that the norm of the operator  $R(\e):{\YY}^2\to    {\YY}^2$  can be estimated as
$\|R(\e)\|<1$ $(\forall \e\in (0,\oo\e] )$.\end{lemma}
\par
{\it Proof}. We have \baaa \left|\sum_{i,j=1}^{n}\w
b_{ij}(x,t) \frac{\p^2 v}{\p x_i\p x_j}(x,t)\right|^2 =
\biggl|\sum_{k\in \N}\biggl( \sum_{i\in \N} \w b_{ki}(x,t)
\frac{\p^2 v}{\p x_k\p x_i}(x,t) +2\sum_{i\notin \N}\w b_{ki}
(x,t) \frac{\p^2 v}{\p x_k\p x_i}(x,t) \biggl)\biggr|^2\nonumber\\
\le \biggl|\sum_{k\in \N}\biggl\{ \biggl( \sum_{i\in \N, i\ne k}
\w b_{ki} (x,t)^2+4\sum_{i\notin \N}\w b_{ki}(x,t)^2
  +\left[1- \frac{\gamma_k}{2}\right]^{-1}\w b_{kk}(x,t)^2\biggr)^{1/2}
  \nonumber\\\times
\biggl( \sum_{i=1,\ldots,n, \ i\ne k} \left|\frac{\p^2 v}{\p x_k\p
x_i}(x,t)\right|^2  +\left[1-\frac{\gamma_k}{2}\right]
\left|\frac{\p^2 v}{\p x_k^2}(x,t)\right|^2 \biggr)^{1/2}
\biggr\}\biggr|^2 \nonumber\\\le \sum_{k\in\N}\biggl(
\sum_{k\in\N, i\ne k} \w b_{ki} (x,t)^2+ 4\sum_{i\notin \N}\w
b_{ki}(x,t)^2 +\left[1- \frac{\gamma_k}{2}\right]^{-1}\w
b_{kk}(x,t)^2 \biggr) \nonumber\\\times
 \sum_{k\in \N}\biggl(
\sum_{i=1,\ldots,n,\,i\ne k} \left|\frac{\p^2 v}{\p x_k\p
x_i}(x,t)\right|^2+\left[1- \frac{\gamma_k}{2}\right]
\left|\frac{\p^2 v}{\p
x_k^2}(x,t)\right|^2\biggr).\,\hphantom{xxx}
\eaaa Hence
$$
\!\left\|\sum_{i,j=1}^n\w b_{ij}\frac{\p^2 v}{\p x_i\p x_j}
\right\|^2_{X^0}\le \w \nu\biggl( \sum_{k\in
\N}\frac{1}{2\gamma_k} \biggr)^{-1}\|v\|^2_{\w X^2}\!<\!
\delta^2\biggl( \sum_{k\in
\N}\!\frac{1}{2\gamma_k}\biggr)^{-1}\!\|v\|^2_{\w X^2}.
$$
In addition, Condition \ref{cond3.1.B} and the embedding theorems
for Sobolev spaces imply the estimates
\begin{eqnarray}
\biggl\|\sum_{i,j=1}^n\bigl( \oo b_{ij}-b^{(\e)}_{ij}\bigr)
\frac{\p^2 v}{\p x_i\p
x_j}\biggr\|_{X^0}+\biggl\|\sum_{i=1}^n\bigl(f_i-f_i^{(\e)}\bigr)
\frac{\p v}{\p x_i}\biggr\|_{X^0}
 +\left\|\bigl(
\lambda-\lambda^{(\e)}\bigr)v\right\|_{X^0}\nonumber\\  \le
C\bigl(\nu_b(\e)+\nu_f(\e)+\nu_{\lambda}(\e)\bigr)\|v\|_{\w X^2} ,
\nonumber\end{eqnarray} where the constant $C$ depends only on
$n$. This proves Lemma \ref{lemma3.3.2}. $\Box$
\par
Let us now complete the proof of Theorem  \ref{Th3.1.2}. By Lemma
\ref{lemma3.3.2}, $(I-R(\e))^{-1}:{\YY}^2\to {\YY}^2$ is a
continuous operator. Let \be
\varphi_{\e}(x,t)\defi \varphi(x,t)e^{K(\e)t}. \label{3.3.15} \ee The
function  $u(x,t)$ is the desired solution of problem
(\ref{3.1.2}), if relation   (\ref{3.3.1}) holds, where \be
u_{\e}=(I-R(\e))^{-1}\bigl[L(\e)\varphi_{\e} +\LL(\e)\Phi \bigr]
\label{3.3.16} \ee because we have
$$
u_{\e}\defi L(\e)\varphi_{\e}+\LL(\e)\Phi +R(\e)u_{\e}.
$$
in view of (\ref{3.3.15})--(\ref{3.3.16}). Therefore,
$$
\|u_{\e}\|_{{\YY}^2}\le
(1-\|R(\e)\|)^{-1}(\|L(\e)\|
\left\|\varphi_{\e}\right\|_{L_2(Q)}+
\|\LL(\e)\|\|\Phi \|_{H^1}) .
$$
This, together with  (\ref{3.3.1}) yields the estimate
(\ref{3.1.6}) and the assertion of Theorem \ref{Th3.1.2}. $\Box$
\section{Uniqueness of a weak solution of It\^o's equation}
Consider the $n$-dimensional vector It\^o's equation \baa
 dy(t)=f(y(t),t)dt+\beta (y(t),t)dw(t),
\\
 y(s)=a.
\label{3.4.1} \eaa By   $y^{a,s}(t)$ we denote a solution of this
equation,  $0\le s\le t\le T$.
\par
In (\ref{3.4.1}), $w(t)$ is a Wiener process of dimension  $n$,
$f(x,t):Q  \to  \R^n$, $\beta (x,t):Q  \to  \R^{n\times n}$,
$Q=\R^n\times (0,T)$ are measurable functions.
\par
Denote
$$     b(x,t)\defi\frac{1}{2} \beta (x,t)\beta (x,t)^{\rm T}.
$$
We assume that the functions   $f(x,t)$, $\beta (x,t)$, $b(x,t)$
are bounded and that the function $b$ satisfies Condition
\ref{cond3.1.A}.
\par
   Let    $(\Omega_0,{\cal F}_0,\P_0)$  be a probability space.
\begin{theorem}\label{Th3.4.1}  (Krylov (1980), Chapter 2). For any random
variable $a \in L^2(\Omega_0,{\cal F}_0,\P_0,\R^n)$, there exists
a set
$$
\biggl\{  (\Omega,{\cal  F},\P),  (w(t),{\cal   F}_t),
y^{a,s}(t) \biggr\},
$$
where    $(\Omega,{\cal F},\P)$ is a probability space  such that
 $a$ $\in $ $L^2 (\Omega,{\cal F},\P)$,
$(w(t),{\cal  F}_t)$  is a Wiener process of dimension  $n$  on
$(\Omega,{\cal F},\P)$,  ${\cal  F}_t\subseteq  {\cal  F}$ is a
filtration of $\sigma $-algebras of events such that  $w(t)-w(s)$
do not depend on $a$ and on ${\cal F}_s$ for  $t>s$, and
$y^{a,s}(t)$ is the solution of (\ref{3.4.1}) for $w(t)$.
\end{theorem}
\par
(In the cited book, the proof was given for non-random $a$, which is
unessential).
\par
  We assume that   $Q=D\times  (0,T)$, where either or $D=\R^n$ or
$D\subseteq \R^n$ is a bounded simply connected domain with
$C^2$-smooth boundary.
\par
Introduce a bounded measurable function $\lambda(x,t):Q \to \C$.
We assume the following condition.
\begin{condition}\label{cond3.4.A}
The functions $b,f,\lambda$ are such that the conclusion of
Theorem  \ref{Th3.1.1} is valid. \end{condition} \begin{remark}
\label{rem3.4.1} It follows from Theorem \ref{Th3.1.1}  that
Condition \ref{cond3.4.A} is satisfied if Condition
\ref{cond3.1.A} is satisfied for $b$, and Condition
\ref{cond3.1.B} is satisfied for $f$ and $\lambda$.
\end{remark}
Let $\chi$ denotes the indicator function.
\begin{theorem}\label{Th3.4.2} Let $a$ be a random vector  with the
probability density function
 $\rho(x)$,  let $a \in  D$  a.s.,  $\rho  \in H^{-1}$,  and $\E|a|^2<  +\infty $.
Let functions   $f(x,t)$, $\beta (x,t)$, $b(x,t)$ be measurable
and bounded, and let Condition \ref{cond3.4.A} be satisfied. Let
$y^{a,s}(t)$ be a weak solution of (\ref{3.4.1}), $
\tau^{a,s}\defi\inf\{t:y^{a,s}(t)\notin D\}$. For the functions $\varphi
\in L_2(Q)$ and $\Phi  \in H^1$, set \baaa
&&F_{a,s}\defi \E\Phi
\bigl(y^{a,s}(T)\bigr)\exp\biggl\{-\int_s^{\tau^{a,s}\land T}
 \lambda(y^{a,s}(r),r)dr\biggr\}\chi_{\{\tau^{a,s}\ge T\}}
\\
&&\hphantom{xxxxxxxxxxxxxxxxxxxx}+\E\int_s^{\tau^{a,s}\land T}\varphi(y^{a,s}(t),t)\exp
\biggl\{-\int_s^t\lambda(y^{a,s}(r),r)dr\biggr\}dt.
\label{3.4.3} \eaaa Then \baaa
 F_{a,s}=\bigl(v(\cdot ,s),\rho\bigr)_{H^0},
 \label{3.4.2'}
\eaaa
where  $v \in Y^2$ is a (unique) solution of problem
(\ref{3.1.2}) for the operator  $A$ given by formula (\ref{3.1.1})
with the above functions $f,b$ and $\lambda$, and
$$
 |F_{a,s} |\le  c \|\rho\|_{H^{-1}}( \|\varphi\|_{L_2(Q)}
+ \|\Phi \|_{H^1}),
$$
where  $c=c({\cal P})$ is a constant occurring in Theorem
\ref{Th3.1.2}.
\end{theorem}
\begin{corollary}\label{corr3.4.1} (The Maximum Principle). Assume
that conditions of Theorem  \ref{Th3.4.2}  are satisfied and, in
addition, that  $\lambda$ is a real function,
 $\varphi(x,t)\ge 0$ for a.e. $x,t$, and $\Phi (x)\ge 0$ for a.e. $
x$. Then the solution  $v$ of problem  (\ref{3.1.2}) is such that  $v(x,t)\ge 0$
for all $t$ for a.e. $x$.
\end{corollary}
\par
    Introduce operators $L_{s,t}:  L_2(D\times
(s,t))  \to H^1$,
${\cal     L}_{s,t}:H^1     \to      H^1$  such that
$v(\cdot ,s)=L_{s,t}\varphi+\LLL_{s,t}\Phi $ is the solution
of the problem
$$
 \left\{\ba
 \frac{\p v}{\p r}(x,r) + Av(x,r)=-\varphi(x,r) , \quad
r<t,
\\
v(x,r)|_{x\in \p D}=0,\quad  v(x,t)=\Phi (x). \ea \right.$$ at the
instant $s$, where $s<t$. By Theorem \ref{Th3.1.1}, these linear
operators are continuous. The conjugate operators
$$L^*_{s,t}:H^{-1} \to L_2(D\times [s,t]),\quad
\LLL^*_{s,t}:H^{-1} \to H^{-1}$$ are also linear and continuous.
\begin{theorem}\label{Th3.4.3}  Under the assumptions of Theorem
\ref{Th3.4.2} (with $D=\R^n$), the weak solution $y^{a,0}(t)$  of
Eqn. (\ref{3.4.1}) with $s=0$  has the probability density
function $p(\cdot,t) \in H^{0}$ for a.e.  $t$. Moreover, $p \in
L_2(Q)$, $p(\cdot,t) \in H^{-1}$ for all  $t$, $p(\cdot ,t)={\cal
L}^*_{0,t}\rho$  and $p=L^*_{0,T}\rho$ for the operators ${\cal
L}^*_{0,t}$, $L^*_{0,T}$ defined for  $\lambda\equiv 0$ (i.e., the
probability density function $p(\cdot,t)$ is uniquely defined as
an element of  $L_2(Q)$  and is uniquely defined as an element of
$H^{-1}$) for all $t$.
\end{theorem}
\par
{\it Proof of Theorems \ref{Th3.4.2}--\ref{Th3.4.3}}. It suffices
to consider
 $s=0$.\par
   (i) Let $\varphi$ and $\Phi $  be such that
$$
               v\defi L\varphi+\LL\Phi  \in C^{2,1}(Q).
$$
Here  $L:X^0 \to {Y}^2$, $\LL:H^1 \to {Y}^2$ are operators such
that   $v=L\varphi+\LL\Phi $ is the solution  of problem \be
\left\{\ba \frac{\p v}{\p t}+Av=-\varphi, \\ v(x,t)|_{x\in\p D}=0,
\quad v(x,T)=\Phi (x) \ea\right.\label{3.4.4} \ee (or the
corresponding Cauchy problem for  $D=\R^n$). In this case relation
(\ref{3.4.2'}) follows from the It\^o formula.
\par
(ii) Let $\varphi  \in  X^0$
and  $\Phi   \in  H^1$ be arbitrary. Introduce the sets $$
        S_1\defi\{\varphi    \in    X^0:L\varphi     \in
C^{2,1}(Q)\},\qquad   S_2\defi\{\Phi  \in H^1: \ \LL\Phi
\in C^{2,1}(Q)\}.
$$
By Theorem \ref{Th3.1.1}, arbitrary functions $\varphi  \in  X^0$
and $\Phi\in  H^1$ can be approximated in these spaces by
$\varphi_{\e}\defi -\p u^{(\e)}/\p t-Au^{(\e)}$ and $\Phi
_{\e}\defi u^{(\e)}(\cdot ,T)$ respectively, where $u^{(\e)}$ is the
Sobolev average of the functions $u=L\varphi$ or $u=\LL\Phi $
respectively: by Theorem \ref{Th3.1.1}, $\varphi_{\e} \to \varphi$
in $X^0$ and $\Phi _{\e} \to  \Phi $ in $H^1$ as $\e \to 0$.
Hence, the sets $S_1$ and  $S_2$ are dense in $X^0$ and in $H^1$,
respectively.
\par
Let $ \oo p\defi L^*_{0,T}\rho$. This is an element of  $X^0$, and
$\oo p(\cdot ,t)={\cal L}^*_{0,t}\rho \in H^{-1} $ for all $t$.
Let  $p(x,t)$ be the probability density function  of the process
$y^{a,0}(t)$ being killed at  $\p D$ if $D\neq\R^n$ and being
killed inside $D$ with the rate $\lambda$. The density   $p(x,t)$
exists by the estimates from  Section 2.3 from Krylov (1980). As
was proved above for $\varphi \in S_1$ and $\Phi  \in S_2$, we
have
$$
       (v(\cdot ,0),\rho)_{H^0}=(\varphi,p)_{X^0}
+(p(\cdot ,T),\Phi )_{H^0}=(\oo  p,\varphi)_{X^0}
+(\oo p(\cdot ,T),\Phi )_{H^0}.
$$
Therefore,
$p=\oo p$ and $p \in X^0$, $p(\cdot ,T)=\oo p(\cdot ,T) \in  H^{-1}$.
\par
Let  $\varphi  \in  X^0$ and   $\Phi    \in   H^1$ be arbitrary,
and let $v\defi L\varphi+\LL\Phi $. Let $v^{(\e)}$  be the Sobolev
average of the function $v$ in $\R^n\times \R$, let
$\varphi_{\e}\defi -\p v^{(\e)} \p t-Av^{(\e)}$, and  let $\Phi
_{\e}\defi  v^{(\e)}(\cdot ,T)$. By Theorem \ref{Th3.1.1}, $\varphi_{\e}
\to \varphi$ in $X^0$ and $\Phi  _{\e}  \to  \Phi  $  in $H^1$ as
$\e \to  0$. We finally obtain the assertion of the theorem from
the relation
$$
\ba
 (v(\cdot ,0),\rho)_{H^0} = \lim_{\e\to 0} (v^{(\e)}(\cdot
,0),\rho)_{H^0} &=\lim_{\e\to 0} ((\varphi_{\e},p)_{X^0}+(p(\cdot
,T),\Phi _{\e})_{H^0})\\& = (p,\varphi)_{X^0}+(p(\cdot ,T),\Phi
)_{H^0}=F_{a,0}. \ea
$$
$\Box$
\begin{theorem}\label{Th3.4.4}
Let  $a$ be a random vector, let $\E|a|^2< +\infty $,  and let
$\rho$ be the probability density function of  $a$, $\rho \in
H^{-1}$. Assume that Condition \ref{cond3.4.A} is satisfied if $f$
is replaced for  $f\equiv 0$, an assume that the function  $f$ is
measurable and bounded. Then problem (\ref{3.4.1}) has a unique
weak solution (i.e., the solution of (\ref{3.4.1}) is univalent
with respect to the probability distribution).
\end{theorem}
\par
{\it Proof}. It suffices to prove the uniqueness of the distribution of the process
$$
z(t)^\top=\bigl[\arctg  y^{a,0}_1 (t),\ldots,\arctg y^{a,0}_n (t)\bigr],
$$
because the function $\arctg :\R  \to  (-\pi,\pi)$  is
one-to-one. We consider  $z(t)$ as a generalized random process defined in Hida (1980)
with the parameter space $L^2([0,T],\oo   {\cal
B}_1,\ell_1,\R^n)$. As is shown in Hida (1980), the
distribution of the process   $z(\cdot)$ is uniquely defined by the
values of the functional
$$
 \ww F_{a,0}(\xi )\defi\E\exp\biggl\{ -\int_0^T i\xi (t)^\top z(t)dt \biggr\}.
$$
on the set  $\xi  \in
L^2([0,T], \oo{\cal B}_1,{\ell}_1, \R^n)$ or on the set of
functions $C([0,T];\R^n)$, which is dense in
$L^2((0,T),\oo{\cal B}_1, {\ell}_1, \R^n)$. Here $i=\sqrt{-1}$.
\par
It is easy to see that
$$
 \ww F_{a,0}(\xi )  = 1-i\E\int_0^T \xi(t)^\top z(t)
\exp\biggl\{ -\int_0^t i \xi (r)^{\rm T}z(r)dr\biggr\} dt.
$$
\par
We first assume that $f\equiv 0$.
By Theorem \ref{Th3.4.2},
$$
\ww F_{a,0}(\xi )=1-i(V,\rho)_{H^0},
$$
where $V=L\varphi$  for
$$
\varphi(x,t)\equiv\xi (t)^\top\bigl[ \arctg  x_1,\ldots,\arctg
x_n\bigr]^\top, \quad \lambda(x,t)\equiv i\varphi(x,t).
$$
Hence    $\ww    F_{a,0}$   is unique for $\xi  \in C((0,T);\R^n)$, and the weak solution is
unique if $f\equiv 0$.\par
 Let $f$ be an arbitrary measurable bounded function. We apply Girsanov theorem.
 Consider the equation
$$
 \left\{\ba   d\ww y(t)=\beta (\ww y(t),t)dw(t),  \\ \ww y(0)=a.
 \ea\right.
$$
As proved above, it has a unique weak solution. By Theorem 2 from  Chapter 3 of
Gihman
and Skorohod (1975), the distribution of the solution $y^{a,0}(t)$  is uniquely
determined by the distribution of $\ww y(t)$.  Hence, the
distribution of  $y^{a,0}(t)$ is defined uniquely. This completes
the proof. $\Box$
\section*{References}
$\hphantom{xxx}$ \par
Anulova, S.V., Krylov, N.V., Liptser, R.S.,
Shyriaev, A.N., and Veretennikov, A.Yu. (1998).
{\it Probability III. Stochastic Calculus.} Berlin: Springer.
\par
Cordes,  H.  O.  (1956). Uber die    erste Randwertaufgabe bei
quasilinerian Differentialgleichungen zweiter Ordnung in mehr als
zwei Variablen. {\it  Math. Ann.} {\bf 131} (iss 3), 278-312.
\par
Dokuchaev, N.G. (1996) On  estimates  for  Ito  processes with
discontinuous disturbances in diffusions. In: {\it Probability
Theory and Mathematical  Statistics (Proc.  Kolmogorov seminars)}.
Gordon and Briach, London, 1996. 133-140.
\par
Gihman, I.I., and Skorohod, A.V. (1975). {\it The Theory of
Stochastic Processes.} Vol. 2. Springer-Verlag, New York.
\par
Hida, T. (1980) {\it Brownian Motion}. Springer-Verlag, 1980.
\par
Kalita E.A. (1989) Regularity of solutions of
Cordes-type elliptic systems of any order. {\it Doklady Acad. Sci.
Ukr. SSR, A}, {\bf  5}, 12-15.
\par Koshelev, A.I. (1982) On exact conditions of regularity of
solutions of for elliptic systems and  Liouville theorem. {\it
Dokl. Akad. Nauk. SSSR}, 265, Iss.6, 1309-1311.
\par
 Krylov, N.V.
(1980).  {\em Controlled Diffusion Processes}.  Shpringer-Verlag.
Berlin, Heidelberg, New York.
\par
 Ladyzhenskaya, O. A., and Ural'tseva, N.N.  (1968). {\it
  Linear and quasilinear elliptic equations}. New York: Academic Press.
\par
Ladyzenskaia, O.A.,  Solonnikov, V.A., and   Ural'ceva, N.N.
(1968). {\it Linear and quasi--linear equations of parabolic
type.} American Mathematical Society. Providence, R.I.
\par
Ladyzhenskaia, O. A. (1985) {\it The boundary value problems of
mathematical physics}. New York: Springer-Verlag.
\par
Landis, E.M. (1998) {\it Second Order equations of elliptic and
parabolic Type}, vol. 171, Amer. Math. Soc., Providence, R.I.,
English transl. in Translations of Math. Monographs.
\par
Liptser, R.S., and Shiryaev, A.N. (2000) {\em Statistics of Random
Processes. I. General Theory}, Berlin, Heidelberg, New York:
Springer-Verlag (2nd ed).
\par Talenti, G. (1965) Sopra una classe di equazioni
elilitticche a coefficienti misurabili. {\it Ann. Math. Pure.
Appl} {\bf 69}, 285-304.
\end{document}